\documentclass[12pt]{article} \usepackage{amssymb,amsmath}

\newtheorem{theorem}{Theorem}[section] \newtheorem{lemma}[theorem]{Lemma}

\newtheorem{corollary}[theorem]{Corollary}

\newenvironment{definition}[1][Definition]{\begin{trivlist} \item[\hskip
\labelsep {\bfseries #1}]}{\end{trivlist}}
\newenvironment{example}[1][Example]{\begin{trivlist} \item[\hskip \labelsep
{\bfseries #1}]}{\end{trivlist}}
 \newcommand{\qed}{\nobreak \ifvmode \relax
\else  \ifdim\lastskip<1.5em \hskip-\lastskip \hskip1.5em plus0em minus0.5em
\fi \nobreak \vrule height0.75em width0.5em depth0.25em\fi}
\newcommand{\E}{\textup{\bf E}}

\begin{document}

\begin{centering}

{\bf \Large A Class of Strongly Homotopy Lie Algebras}

{\bf \Large \vspace{.5cm}with Simplified sh-Lie Structures}
\\

\renewcommand{\thefootnote}{\fnsymbol{footnote}}
\vspace{1.4cm}
{\large Samer Al-Ashhab}\\

\vspace{.5cm}

Department of Mathematics, North Carolina State University, Raleigh,
NC 27695-8205.\\
\it E-mail: ssalash@unity.ncsu.edu\\
\rm
\vspace{.5cm}

\rm
\vspace{.5cm}

\begin{abstract}

It is known that a single mapping defined on one term of a differential graded 
vector space extends to a strongly homotopy Lie algebra structure on the graded
space when that mapping satisfies two conditions. This strongly homotopy Lie
algebra is nontrivial (it is not a Lie algebra); however we show that one can
obtain an sh-Lie algebra where the only nonzero mappings defining it are the
lower order mappings. This structure applies to a significant class of
examples. Moreover in this case the graded space can be replaced by another
graded space, with only three nonzero terms, on which the same sh-Lie structure
exists.

\end{abstract}

\end{centering}

\section{Introduction and sh-Lie algebras} Strongly homotopy Lie algebras
(sh-Lie algebras/structures) have recently been the focus of study in
mathematics \cite{BFLS98,LM95,LS93}. Their applications have appeared in
mathematics \cite{RW98}, in mathematical physics \cite{BFLS98,FLS01}, and in
physics \cite{Z93}. We first present some background material. Then in section
2 we prove the main result. Finally in section 3, we present three examples,
two of which have appeared in two different applications.

To begin with our discussion, let $X_*$ be a graded vector space with a
differential (lowering the degree by 1) $l_1$, and maps $\eta:X_n \to H_n,
\lambda:H_n \to X_n$, and $s:X_n \to X_{n+1}$ (i.e. $s$ raises the degree by
1), where $H_*$ is the homology complex of the complex $X_*,  H_0$ is generally
non-trivial and $H_n=0$ for $n>0.$ We also assume that \begin{equation}
\label{s} \lambda \circ \eta - 1_{X_*} = l_1\circ s + s\circ l_1.
\end{equation} Notice that $\ref {s}$ simplifies to $\lambda \circ \eta -
1_{X_*} = l_1\circ s$ in degree 0, and $-1_{X_*} = l_1\circ s + s\circ l_1$ in
degree $>0$. For more details we refer the reader to \cite{BFLS98}. The
following is the formal definition of sh-Lie structures (see \cite{BFLS98})

\begin{definition} An sh-Lie structure on a graded vector space $X_*$ is a
collection of linear, skew-symmetric maps $l_k: \bigotimes ^k X_* \to X_*$ of
degree $k-2$ that satisfy the relation $$\sum_{i+j=n+1}\sum_{unsh(i,n-i)}
e(\sigma)(-1)^ \sigma (-1)^{i(j-1)}l_j(l_i
(x_{\sigma(1)},\cdots,x_{\sigma(i)}),\cdots,x_{\sigma(n)}) = 0,$$ where $ 1
\leq i,j$. \end{definition} Notice that in this definition $e(\sigma)$ is the
Koszul sign which depends on the permutation $\sigma$ as well as on the degree
of the elements $x_1,x_2,\cdots,x_n$ (where a minus sign is introduced whenever
two consecutive odd elements are permuted, see \cite{LM95}). Also observe that
it is convenient in this context to suppress some of the notation and assume
the summands are over the appropriate unshuffles with their corresponding signs
(e.g. if $n=3$ one writes $l_1l_3 + l_2l_2 + l_3l_1 = 0$). \\ 

We assume the existence of a linear skew-symmetric map $\tilde{l}_2: X_0
\bigotimes X_0 \to X_0$ satisfying conditions $(i)$ and $(ii)$ below so that an
sh-Lie structure exists, where we quote the following from \cite{BFLS98}

\begin{theorem} \label{shLie}
A skew-symmetric linear map $\tilde{l}_2: X_0 \bigotimes X_0 \to
X_0$ that satisfies conditions (i) and (ii) below extends to an sh-Lie
structure on the graded space $X_*$;\\ $(i) \quad \tilde{l}_2(c,b_1) = b_2$ \\
$\displaystyle (ii) \sum _{\sigma \in unsh(2,1)} (-1)^\sigma
\tilde{l}_2(\tilde{l}_2(c_{\sigma(1)},c_{\sigma(2)}),c_{\sigma(3)}) =b_3$ \\
where $c, c_1, c_2, c_3$ are cycles and $b_1, b_2, b_3$ are boundaries in
$X_0$. \end{theorem}

\section{The main theorem}
While Theorem $\ref {shLie}$ guarantees the existence of an sh-Lie structure on
the graded vector space $X_*$, we \textit{show} that one can always choose an
sh-Lie structure such
that \\ 1. $l_2=0$ in degree $>1$.\\ 2. $l_3=0$ in degree $>0$.\\
3. $l_n \equiv 0,n>3.$

\bigskip

\noindent{\bf Remark} Markl has observed (see \cite{BFLS98}) that such an sh-Lie
structure exists in the case that $l_2(c,b)=0$ for each cycle $c$ and boundary
$b.$ To our knowledge even the proof of this special case has not been
published although the result is known to specialists in the area. \bigskip

We find it convenient in what follows to refer to the image of a combination of
maps by the combination itself: for example $l_2l_1$ would stand for the image
of the map $l_2l_1$
acting on some element in the appropriate space as the context implies (e.g
$l_2l_1$ may stand for $l_2l_1(x_p \otimes x_q)= l_2(l_1x_p \otimes x_q +
(-1)^p x_p \otimes l_1x_q) = l_2(l_1x_p \otimes x_q) + (-1)^pl_2(x_p \otimes
l_1x_q)$ for $x_p \in X_p$ and $x_q \in X_q$). $l_2l_2$ would stand for the
three unshuffles of the composite $l_2l_2$, again acting on an element in the
approprite space (where we skip writing the element it's acting on, e.g.
$l_2l_2= l_2l_2(x_0\otimes x_0'\otimes x_0'')= l_2(l_2(x_0,x_0'),x_0'')-
l_2(l_2(x_0,x_0''),x_0')+l_2(l_2(x_0',x_0''),x_0)$ for $x_0\otimes x_0'\otimes
x_0'' \in X_0\bigotimes X_0\bigotimes X_0$) ...etc. Let's also quote the
following from \cite{BFLS98} as it is needed in our proof:
\begin{lemma} (i) $l_2l_1$ is a boundary.\\ (ii) $l_2l_2 +
l_3l_1$ is a boundary.\\ (iii) More generally $\displaystyle
(\sum_{i+j=n+1,j>1} (-1)^{i(j-1)}l_jl_i)$ is a boundary. \end{lemma}

To begin with the proof we define $l_2$ inductively by $$l_2=-s\circ l_2 l_1$$
(this is just $-sl_2l_1$, i.e. $s$ acts on the image of $l_2l_1$) where we
begin with $l_2=\tilde{l}_2$ in degree 0, and recall that $s$ satisfies
$$\lambda \circ \eta - 1_{X_*} = l_1\circ s$$ in degree 0, and $$-1_{X_*} =
l_1\circ s + s\circ l_1$$ in degree $>0$ (see $\ref {s}$). One checks that if
$l_2l_1$ is in degree 0 then $$l_1(-sl_2l_1)=l_2l_1- \lambda \circ
\eta(l_2l_1)=l_2l_1,$$ where $(\lambda \circ \eta)(l_2l_1)=0$ since $l_2l_1$ is
a boundary. While if $l_2l_1$ is in degree $>0$ then $$l_1(-sl_2l_1) = l_2l_1 +
sl_1(l_2l_1)= l_2l_1,$$ where $l_1(l_2l_1)=0$ since $l_2l_1$ is a boundary. So
we have a well-defined chain map $l_2$ satisfying $l_1l_2=l_2l_1$. Now we show
that $l_2$ as defined above is {\em zero} in degree $>1$. 

First consider $l_2$ on $\underline {X_1\otimes X_1}$: take an element
$x_1\otimes x_1' \in X_1\otimes X_1$. We have $l_2(x_1\otimes x_1')=
-s\{l_2(l_1x_1\otimes x_1' - x_1\otimes l_1x_1')\}.$ But $l_2(l_1x_1\otimes
x_1')=l_2(x_1\otimes l_1x_1')$, since by definition $$l_2(l_1x_1\otimes x_1')=
-sl_2l_1(l_1x_1\otimes x_1')=-sl_2(l_1x_1\otimes l_1x_1'),$$ and
$$l_2(x_1\otimes l_1x_1')= -sl_2l_1(x_1\otimes l_1x_1')=-sl_2(l_1x_1\otimes
l_1x_1').$$ So $l_2=0$ on $X_1\otimes X_1$. Now consider $l_2$ on $\underline
{X_2\otimes X_0}$: take $x_2\otimes x_0 \in X_2\otimes X_0$. Then
$l_2(x_2\otimes x_0)=-sl_2(l_1x_2\otimes x_0)$, but $l_2(l_1x_2\otimes
x_0)=-sl_2l_1(l_1x_2\otimes x_0)=-sl_2(0)=0$. So $l_2=0$ on $X_2\otimes X_0$.

Proceeding by induction one then shows that $l_2=0$ on $X_n\otimes X_0$ with $n
\geq 3$: $l_2(x_n\otimes x_0)=-sl_2l_1(x_n\otimes x_0) =-sl_2(l_1x_n\otimes
x_0)=-s(0)=0$. On the other hand consider $l_2$ on $X_n\otimes X_m, n>1, m \geq
1$: $l_2(x_n\otimes x_m)=-sl_2l_1(x_n\otimes x_m) =-sl_2(l_1x_n\otimes x_m
+(-1)^n x_n\otimes l_1x_m)=-s(0)=0$. This way one has $l_2 \equiv 0$ in degree
$> 1$.
Now consider the map $l_3$. Define $$l_3=s\circ l_2l_2,$$ in degree 0 and
then inductively by $$l_3=s\circ (l_2l_2+l_3l_1),$$ in degree $>0$. One checks
that in degree 0 $$-l_1l_3=-l_1s(l_2l_2)=l_2l_2- (\lambda \circ
\eta)(l_2l_2)=l_2l_2,$$ where $(\lambda \circ \eta)(l_2l_2)=0$ since $l_2l_2$
is a boundary. While in degree $>0$ we have $$-l_1l_3=-l_1(s(l_2l_2+l_3l_1) =
l_2l_2+l_3l_1 + s\{l_1(l_2l_2+l_3l_1)\}= l_2l_2+l_3l_1,$$ where
$l_1(l_2l_2+l_3l_1)=0$ since $l_2l_2+l_3l_1$ is a boundary. So we have a well
defined chain map $l_3$ satisfying $l_1l_3+l_3l_1+l_2l_2=0$.

Now consider $l_3$ on $X_1\otimes X_0\otimes X_0$: take an element $x_1\otimes
x_0\otimes x_0' \in X_1\otimes X_0\otimes X_0$. By definition we have: \\
$(l_2l_2+l_3l_1)(x_1\otimes x_0\otimes x_0')=$\\ $-sl_2l_1(l_2(x_1,x_0),x_0')+
sl_2l_1(l_2(l_2(x_1,x_0'),x_0)-sl_2l_1(l_2(l_2(x_0,x_0'),x_1)$\\
$+s\{l_2(l_2(l_1x_1,x_0),x_0')-
l_2(l_2(l_1x_1,x_0'),x_0)+l_2(l_2((x_0,x_0'),l_1x_1)\}=$\\ \\
$-sl_2(l_2(l_1x_1,x_0),x_0')+
sl_2(l_2(l_2(l_1x_1,x_0'),x_0)-sl_2(l_2(l_2(x_0,x_0'),l_1x_1)$\\
$+sl_2(l_2(l_1x_1,x_0),x_0')-
sl_2(l_2(l_1x_1,x_0'),x_0)+sl_2(l_2((x_0,x_0'),l_1x_1)=0$.\\ \\ So
$l_3=s\{(l_2l_2+l_3l_1)\}=s(0)=0$ on $X_1\otimes X_0\otimes X_0$. 

One can then proceed by induction to find that $l_3=0$ in degree $>1$, for
example take $x_1\otimes x_1'\otimes x_0$ in $X_1\otimes X_1\otimes X_0$.
$l_3(x_1\otimes x_1'\otimes x_0)=s(l_2l_2+l_3l_1)(x_1\otimes x_1'\otimes x_0)=
s\{l_2(l_2(x_1,x_1'),x_0)-l_2(l_2(x_1,x_0),x_1')+l_2(l_2(x_1',x_0),x_1)+
l_3(l_1x_1\otimes x_1'\otimes x_0)+l_3(x_1\otimes l_1x_1'\otimes x_0)\}=0$,
since $l_3$ and $l_2$ are {\em zero} in degrees 1 and 2 respectively.
Now consider the map $l_4$. Define $$l_4=s\circ (l_3l_2-l_2l_3),$$ in degree
0, and then inductively by $$l_4=s\circ (l_3l_2-l_2l_3-l_4l_1),$$ in degree
$>0$. As before one can easily check that $l_4$ is a well-defined map that
satisfies the corresponding sh-Lie relation at this step of the construction
(i.e. $l_1l_4 - l_4l_1 + l_3l_2 - l_2l_3 = 0$).

Consider the value of $l_4$ on $x_0\otimes x_0'\otimes x_0''\otimes x_0''' \in
X_0\otimes X_0\otimes X_0\otimes X_0$. We find it convenient in the following
calculation to use the identity $l_3(l_2(y_0,y_0'),y_0'',y_0''')=
(sl_2l_2)(l_2(y_0,y_0'),y_0'',y_0''')= s\{l_2(l_2(l_2(y_0,y_0'),y_0''),y_0''')-
l_2(l_2(l_2(y_0,y_0'),y_0'''),y_0'')+l_2(l_2(y_0'',y_0'''),l_2(y_0,y_0'))\}$,
for $y_0,y_0',y_0'',y_0''' \in X_0$. By definition $l_4$ is the values of $s$
on: \\ $l_3(l_2(x_0,x_0'),x_0'',x_0''')-l_3(l_2(x_0,x_0''),x_0',x_0''')+
l_3(l_2(x_0,x_0'''),x_0',x_0'')+ \\ l_3(l_2(x_0',x_0''),x_0,x_0''')-
l_3(l_2(x_0',x_0'''),x_0,x_0'')+l_3(l_2(x_0'',x_0'''),x_0,x_0')- \\
l_2(l_3(x_0,x_0',x_0''),x_0''')+l_2(l_3(x_0,x_0',x_0'''),x_0'')-
l_2(l_3(x_0,x_0'',x_0'''),x_0')+ \\ l_2(l_3(x_0',x_0'',x_0'''),x_0)=$\\ 

\noindent $s\{l_2(l_2(l_2(x_0,x_0'),x_0''),x_0''')-
l_2(l_2(l_2(x_0,x_0'),x_0'''),x_0'')+l_2(l_2(x_0'',x_0'''),l_2(x_0,x_0'))\}-
s\{l_2(l_2(l_2(x_0,x_0''),x_0'),x_0''')- l_2(l_2(l_2(x_0,x_0''),x_0'''),x_0')+
l_2(l_2(x_0',x_0'''),l_2(x_0,x_0''))\} +
s\{l_2(l_2(l_2(x_0,x_0'''),x_0'),x_0'')-
l_2(l_2(l_2(x_0,x_0'''),x_0''),x_0')+l_2(l_2(x_0',x_0''),l_2(x_0,x_0'''))\}+
s\{l_2(l_2(l_2(x_0',x_0''),x_0),x_0''')-
l_2(l_2(l_2(x_0',x_0''),x_0'''),x_0)+l_2(l_2(x_0,x_0'''),l_2(x_0',x_0''))\}-
s\{l_2(l_2(l_2(x_0',x_0'''),x_0),x_0'')- l_2(l_2(l_2(x_0',x_0'''),x_0''),x_0)+
l_2(l_2(x_0,x_0''),l_2(x_0',x_0'''))\} +
s\{l_2(l_2(l_2(x_0'',x_0'''),x_0),x_0')-
l_2(l_2(l_2(x_0'',x_0'''),x_0'),x_0)+l_2(l_2(x_0,x_0'),l_2(x_0'',x_0'''))\}+
s\{l_2(l_1l_3(x_0,x_0',x_0''),x_0''')-l_2(l_1l_3(x_0,x_0',x_0'''),x_0'')+
l_2(l_1l_3(x_0,x_0'',x_0'''),x_0')-$\\ $l_2(l_1l_3(x_0',x_0'',x_0'''),x_0)\}=$
\\ \\ \noindent $sl_2\{(l_2l_2+l_1l_3)(x_0,x_0',x_0''),x_0'''\}-
sl_2\{(l_2l_2+l_1l_3)(x_0,x_0',x_0'''),x_0''\}+$\\
$sl_2\{(l_2l_2+l_1l_3)(x_0,x_0'',x_0'''),x_0'\}-
sl_2\{(l_2l_2+l_1l_3)(x_0',x_0'',x_0'''),x_0\}+$\\
$l_2(l_2(x_0'',x_0'''),l_2(x_0,x_0'))+l_2(l_2(x_0,x_0'),l_2(x_0'',x_0'''))-
l_2(l_2(x_0',x_0'''),l_2(x_0,x_0''))-$\\ $l_2(l_2(x_0,x_0''),l_2(x_0',x_0'''))+
l_2(l_2(x_0',x_0''),l_2(x_0,x_0'''))+l_2(l_2(x_0,x_0'''),l_2(x_0',x_0''))=0$
since $l_2l_2+l_1l_3=0$ in degree 0 and $l_2$ is skew-symmetric. So we have
$l_4=0$ in degree 0.

Further $l_4$ is inductively found to be {\em zero} in higher degrees since
$l_2=0$ and $l_3=0$ in degrees $>1$ and $>0$ respectively.\\ 

Next we inductively define for $n>4$, $$l_n=s\circ (\sum_{i,j>1}
(-1)^{i(j-1)}l_jl_i),$$ in degree 0, and $$l_n=s\circ (\sum_{i+j=n+1,j>1}
(-1)^{i(j-1)}l_jl_i),$$ in degree $>0$. (Again these are well defined maps for
the sh-Lie structure.) The combination of maps (the $l_k$'s within the $s$) in
degree 0, and then inductively in degree $>0$, leads to 0 so that one has $l_n
\equiv 0$ for $n>4$ (Notice that for $l_5$ in degree 0 one encounters $l_3l_3$,
the inside $l_3$ raises the degree from 0 to 1 so that the combination is 0).\\
\\  Summarizing: \begin{theorem} \label{simp} Given a graded space $X_*$ and a
skew-symmetric linear map $\tilde{l}_2: X_0 \bigotimes X_0 \to X_0$ that
satisfies conditions (i) and (ii), there exists an sh-Lie structure on $X_*$
such that \\ 1. $l_2=0$ in degree $>1$.\\ 2. $l_3=0$ in degree $>0$.\\ 3. $l_n
\equiv 0,n>3.$ \end{theorem}

\begin{corollary} \label{simpc} Under the same hypotheses in the theorem there
exists an sh-Lie structure on the graded space $$\cdots \longrightarrow 0
\longrightarrow 0 \longrightarrow X_2 \xrightarrow{i} X_1 \xrightarrow{l_1}
X_0,$$ where $X_1, X_0$, $l_1:X_1\rightarrow X_0$, and $l_k, k>1$ are as above,
but with $X_2=ker l_1$ and the inclusion $i:X_2 \rightarrow X_1$.
\end{corollary}

\section{Examples} We will primarily consider three examples. The first of
which fits perfectly into our discussion. It first appeared in \cite{RW98} as
the authors determined an sh-Lie structure on a Courant algebroid in the sense
of the example given below. For convenience we recall the definition of a
Courant algebroid. \begin{definition} A Courant algebroid is a vector bundle $E
\to M$ equipped with a nondegenerate symmetric bilinear form $<\cdot,\cdot>$ on
the bundle, a skew-symmetric bracket $[\cdot,\cdot]$ on $\Gamma (E)$, and a
bundle map $\rho:E \to TM$ such that the following properties are satisfied: \\
1. For any $e_1,e_2,e_3 \in \Gamma (E), J(e_1,e_2,e_3)= \mathcal D
T(e_1,e_2,e_3)$; \\ 2. for any $e_1,e_2 \in \Gamma (E), \rho[e_1,e_2]=[\rho
e_1,\rho e_2]$; \\ 3. for any $e_1,e_2 \in \Gamma (E)$ and $f \in C^\infty(M),
[e_1,fe_2] = f[e_1,e_2] + (\rho(e_1)f)e_2 - <e_1,e_2> \mathcal D f$; \\ 4.
$\rho \circ \mathcal D = 0$; \\ 5. for any $e,h_1,h_2 \in \Gamma (E),
\rho(e)<h_1,h_2>=<[e,h_1]+\mathcal D <e,h_1>,h_2>+<h_1,[e,h_2]+\mathcal D
<e,h_2>>;$ \\ where $$J(e_1,e_2,e_3)=[[e_1,e_2],e_3] + [[e_2,e_3],e_1]
+[[e_3,e_1],e_2],$$ and $T(e_1,e_2,e_3)$ is the function on the base space $M$
defined by $$ T(e_1,e_2,e_3)=\frac{1}{3}<[e_1,e_2],e_3> + c.p.$$ ($c.p.$ here
denotes the cyclic permutations of the $e_i$'s) and $\mathcal D:C^\infty(M) \to
\Gamma (E)$ is defined such that the following identity holds $$<\mathcal D
f,e> = \frac{1}{2} \rho(e) f.$$ \end{definition} For more details on Courant
algebroids we refer the reader to \cite{RW98} and the references therein.
\begin{example} Let $E$ be a Courant algebroid over a manifold $M$, and
consider the sequence $$\cdots \rightarrow 0 \rightarrow \text{ker}
\mathcal D \xrightarrow{i} C^\infty (M) \xrightarrow{\mathcal D} \Gamma (E)$$
where $i:\text{ker} \mathcal D \to C^\infty (M)$ is the inclusion, and one
assumes that $X_0=\Gamma (E), X_1=C^\infty (M)$, and $X_2=\text{ker} \mathcal
D$. Define $\tilde{l}_2(e_1,e_2)=[e_1,e_2]$ (this is just $l_2$ in degree 0).
It was shown in \cite{RW98}
that $\tilde{l}_2$ satisfies condition $(i)$ as in Theorem
$\ref{shLie}$, whereas condition $(ii)$ from the same Theorem follows from the
the first axiom in the above definition yielding an sh-Lie structure. The
authors \cite{RW98}
have also shown that the sh-Lie structure has the explicit formulas
\[ \begin{array}{lccl} l_2(e_1,e_2) & = & [e_1,e_2] & \text{in degree 0} \\
l_2(e,f) & = & <e,\mathcal D f> & \text{in degree 1} \\ l_2 & = & 0 & \text{in
degree} >1 \\ l_3(e_1,e_2,e_3) & = & -T(e_1,e_2,e_3) & \text{in degree 0} \\
l_3 & = & 0 & \text{in degree} > 0 \\ l_n & = & 0 & \text{for } n>3.
\end{array} \] What is interesting about this example is that the only nonzero
maps are the same ones as in Theorem $\ref{simp}$ and Corollary $\ref{simpc}$.
In addition notice that the structure of this complex is similar to the
``simplified" complex that appears in Corollary $\ref{simpc}$. \qed
\end{example}

Our next example comes from Lagrangian field theory, in particular it relates
to the Poisson brackets of local functionals where the sh-Lie structure exists
on a ``DeRahm complex" as in \cite{BFLS98}. We refer the reader to \cite{BFLS98}
and the
references therein for more details regarding this subject. \begin{example} let
$E \to M$ be a vector bundle where the base space $M$ is an $n$-dimensional
manifold and let $J^\infty E$ be the infinite jet bundle of $E.$ Consider the
complex $$\Omega ^{0,0}(J^\infty E) \to \Omega ^{1,0}(J^\infty E) \to \cdots
\to \Omega ^{n-1,0}(J^\infty E) \to \Omega ^{n,0}(J^\infty E)$$ with a
differential $d_H$ which in local coordinates takes the form $d_H = dx^i D_i$,
i.e. if $\alpha = \alpha _I dx^I$ then $d_H \alpha = D_i \alpha _I dx^i \wedge
dx^I$. Here $D_i$ is the total derivative derivation defined on the algebra of
local functions on $J^\infty E. $ It is defined by $\displaystyle D_i =
\frac{\partial}{\partial x^i} + u^a_{iJ}\frac{\partial} {\partial u^a_J}$
(we assume the summation convention, i.e., the sum is over all $a$ and
multi-index $J$). In this case $\tilde{l}_2$ was defined in \cite{BFLS98} by
\begin{equation} \label{l2} \tilde{l}_2(P\nu ,Q\nu)= \omega(\E(P),\E(Q)) \nu =
\omega ^{ab}(\E_b(Q))\E_a(P) \nu \end{equation} where $P\nu ,Q\nu \in \Omega
^{n,0}(J^\infty E)$ and $\E$ is the Euler-Lagrange operator with components
$$\E_a(P)=(-D)_I(\frac{\partial P} {\partial u^a_I}).$$ The bilinear mapping 
$\omega$ is a skew-symmetric total differential operator with the $\omega
^{ab}$'s as its components (see \cite{BFLS98} for more details). It was shown in
\cite{BFLS98} that
$\tilde{l}_2$ satisfies conditions $(i)$ and $(ii)$ as in Theorem $\ref
{shLie}$. In fact $\tilde{l}_2$ satisfies condition $(i)$ in a strong sense
(with 0 on the right-hand side of the equation). Markl noted in \cite{BFLS98}
that with this
strong condition the higher order maps can be chosen to be zero (The result in
this paper's Theorem $\ref{simp}$ is yet stronger since it does not require
that the right-hand side of $(i)$ be zero, only that it be a boundary). Here is
a summary of the structure \[ \begin{array}{lccl} l_2(P\nu,Q\nu) & = &
\omega(\E(Q),\E(P)) \nu & \text{in degree 0} \\ l_2 & = & 0 & \text{in degree}
>0 \\ l_3(P\nu,Q\nu,R\nu) & is & nonzero & \text{in degree 0} \\ l_3 & = & 0 &
\text{in degree} > 0 \\ l_n & = & 0 & \text{for } n>3. \qed \end{array} \]
\end{example}

Our last example is within the context of symplectic manifolds, where we refer
the reader to \cite{AM94}. We include details in this example on how the sh-Lie
structure maps are obtained. Notice that in this example the strong version of
$(i)$ of Theorem $\ref {shLie}$ does not hold in general, but our weaker
hypothesis does hold. \begin{example} Consider the following sequence
$$0\rightarrow R \rightarrow \Omega ^0(M) \rightarrow \Omega_C^1(M) \rightarrow
0,$$ where $\Omega ^0(M)$ is the set of smooth real-valued functions on the
symplectic manifold $(M,\omega)$ and $\Omega_C^1(M)$ is the set of closed one
forms on $M$. We take $X_0=\Omega_C^1(M), X_1=\Omega ^0(M)$ and $X_2=R$. The
chain map is $l_1=i:R \rightarrow \Omega ^0(M)$, and $l_1=d:\Omega ^0(M)
\rightarrow \Omega_C^1(M)$, where $i$ is the inclusion and $d$ is the
differential operator. We then define a bilinear skew-symmetric map
$\tilde{l}_2$ on $X_0 \times X_0$ by $$\tilde{l}_2(\alpha,\beta)=\{\alpha,\beta
\},$$ where \{.,.\} is a Poisson bracket on $\Omega^1(M)$ (e.g. see definition
3.3.7 in \cite{AM94}). Notice that $\tilde{l}_2$ satisfies the two conditions
needed to guarantee the existence of an sh-Lie algebra (Theorem $\ref{shLie}$).

Now to extend $\tilde{l}_2$, first take an element in $X_1 \otimes X_0$ say
$f\otimes \beta$, then $\tilde{l}_2l_1(f\otimes \beta) = \tilde{l}_2(df \otimes
\beta)= \{df,\beta \}$. Now notice that $h=L_{\beta^\#}f + c =
-i_{X_f}\beta+L_{\beta^\#}f+i_{X_f}i_{\beta^\#}\omega +c \in X_1$ satisfies
\begin{eqnarray*} l_1(h) & = & d(-i_{X_f}\beta+L_{\beta^\#}f +
i_{X_f}i_{\beta^\#}\omega +c) \\ & = &
-di_{X_f}\beta-i_{X_f}d\beta+dL_{\beta^\#}f+d(i_{X_f}i_{\beta^\#}\omega)\\ & =
& -L_{X_f}\beta+L_{\beta^\#}df+d(i_{X_f}i_{\beta^\#}\omega)\\ & = & \{df,\beta
\}, \end{eqnarray*} where $\omega$ is the {\em symplectic} 2-form on $M$. So if
we take $c=0$ we get $l_2(f\otimes\beta) = L_{\beta^\#}f$. Then $l_2$
would be defined on $X_0\otimes X_1$ by skew-symmetry. 
To proceed take an element in $X_1\otimes X_1$ say $f\otimes g$, and notice that
$l_2l_1(f\otimes g) =
l_2(df \otimes g - f\otimes dg)=
L_{df^\#}g - (-L_{dg^\#}f) = L_{X_f}g + L_{X_g}f = 0$. Hence $l_2$ on
$X_1\otimes X_1$ is {\em zero}.
Now take an element in $X_2\otimes X_0$ say $k\otimes \beta$, then
$l_2l_1(k\otimes \beta) = l_2(k\otimes \beta + 0)= L_{\beta^\#}k = 0$, since
$k$ is a constant function. Therefore $l_2$ is {\em zero} on
$X_2\otimes X_0$. By skew-symmetry $l_2$ will also be {\em zero} on
$X_0\otimes X_2$.

Next we turn to $l_3$, take an element in $X_0\otimes X_0 \otimes X_0$, say
$\alpha \otimes \beta \otimes \gamma$, then $l_2l_2$ maps it into
$\{\{\alpha,\beta\},\gamma\} -\{\{\alpha,\gamma\},\beta\} +
\{\{\beta,\gamma\},\alpha\})$ which is 0 (The Jacobi identity). So we have
$l_3=0$ on $X_0\otimes X_0 \otimes X_0$. Now take $f\otimes \beta \otimes
\gamma \in X_1\otimes X_0 \otimes X_0$. Under $l_2l_2 + l_3l_1$ it is mapped to
$$l_2(L_{\beta^\#}f,\gamma) - l_2(L_{\gamma^\#}f,\beta) +
l_2(\{\beta,\gamma\},f) + 0 =$$ $$L_{\gamma^\#}L_{\beta^\#}f -
L_{\beta^\#}L_{\gamma^\#}f - L_{\{\beta,\gamma\}^\#}f=$$
$$L_{\{\beta,\gamma\}^\#}f-L_{\{\beta,\gamma\}^\#}f=0.$$ So $l_3$ is {\em zero}
on $X_1 \otimes X_0 \otimes X_0$. Utilizing skew-symmetry we have $l_3 = 0$ in
degree 1. On higher degrees $l_3$ is trivially {\em zero} since $X_n=0$ for
$n>2$, and furthermore $l_n=0$ for $n>3$. To summarize \[ \begin{array}{lccl}
l_2((\alpha,\beta) & = & \{\alpha,\beta \} & \text{in degree 0} \\ l_2(f,\beta)
& = & L_{\beta^\#}f & \text{in degree 1} \\ l_2 & = & 0 & \text{in degree} >1
\\ l_n & = & 0 & \text{for } n>2.
\end{array} \] \end{example}

\noindent {\bf Acknowledgements} I would like to thank Professor Ron Fulp for
suggestions, advice, and help as this paper is a part of my PhD research which
he advises. I would also like to thank Professor Tom Lada for discussions about
sh-Lie algebras, in particular discussions about the main result in this paper.

\end{document}